\documentclass[runningheads]{llncs}
\usepackage{graphicx}
\usepackage{amsmath}
\usepackage{amssymb}
\usepackage{booktabs} 
\usepackage{multirow} 
\usepackage{array}
\usepackage{todonotes}
\usepackage{adjustbox} 

\usepackage[hidelinks]{hyperref}

\newcommand\calS{\mathcal{S}}
\newcommand\GDA{\mathrm{GD(\mathcal{A}_p)}}
\newcommand\DHV{\mathrm{\Delta HV(\mathcal{A}_p)}}

\newcommand{\dd}[2]{\frac{\partial {#1}}{\partial {#2}} }

\begin{document}

\title{Multi-objective Optimization by Uncrowded Hypervolume Gradient Ascent}
\titlerunning{Multi-objective Optimization by Uncrowded Hypervolume Gradient Ascent}

%
\author{Timo M. Deist\inst{1}\orcidID{0000-0003-0057-1535}\thanks{These authors contributed equally.} \and
Stefanus C. Maree\inst{1}\textsuperscript{\thefootnote} \and
Tanja Alderliesten\inst{2}
\and
Peter A.N. Bosman\inst{1}}
\authorrunning{T.M. Deist, S.C. Maree et al.}
%

\institute{Centrum Wiskunde \& Informatica, Life Sciences and Health Research Group, Amsterdam, The Netherlands. 
\email{\{timo.deist, maree, peter.bosman\}@cwi.nl}
\and
Leiden University Medical Center, Department of Radiation Oncology, \\ Leiden, The~Netherlands. 
\email{t.alderliesten@lumc.nl}}

\maketitle

\begin{abstract}

Evolutionary algorithms (EAs) are the preferred method for solving black-box multi-objective optimization problems, but when gradients of the objective functions are available, it is not straightforward to exploit these efficiently. 
By contrast, gradient-based optimization is well-established for single-objective optimization. A single-objective reformulation of the multi-objective problem could therefore offer a solution. Of particular interest to this end is the recently introduced uncrowded hypervolume (UHV) indicator, which takes into account dominated solutions.
In this work, we show that the gradient of the UHV can often be computed, which allows for a direct application of gradient ascent algorithms.
We compare this new approach with two EAs for UHV optimization as well as with one gradient-based algorithm for optimizing the well-established hypervolume. On several bi-objective benchmarks, we find that gradient-based algorithms outperform the tested EAs by obtaining a better hypervolume with fewer evaluations whenever exact gradients of the multiple objective functions are available and in case of small evaluation budgets. For larger budgets, however, EAs perform similarly or better. 
We further find that, when finite differences are used to approximate the gradients of the multiple objectives, our new gradient-based algorithm is still competitive with EAs in most considered benchmarks. Implementations are available at \url{https://github.com/scmaree/uncrowded-hypervolume}.

\end{abstract}

\begin{tikzpicture}[remember picture, overlay]
\node[anchor = south, yshift = 10pt] at (current
page.south){\fbox{\parbox{\textwidth}{\centering The final authenticated version is available online at\\ \url{https://doi.org/10.1007/978-3-030-58115-2_13}}}};
\end{tikzpicture}

\section{Introduction}
Evolutionary algorithms (EAs) are the preferred method for solving black-box multi-objective (MO) optimization problems, when assuming the underlying details of the problem are unknown \cite{deb2001book}. However, when gradient information of the objective functions is available, it is not straightforward to exploit this information efficiently in the optimization process. 
This can be mainly attributed to the two-sided goal of multi-objective optimization, which is to obtain a set of solutions, known as an approximation set, on the one hand containing solutions that are (near) Pareto optimal, and on the other hand representing a diverse set of trade-offs between the objectives \cite{bosman2003}. 

When considering a to-be-minimized bi-objective function $\mathbf{f} : \mathbb{R}^n \rightarrow \mathbb{R}^2$, 
the Karush-Kuhn-Tucker (KKT) \cite{KKT51,peitz2017} conditions can be used to identify a descent direction $d(\mathbf{x})$ for a solution $\mathbf{x}\in\mathbb{R}^n$ for which all objectives are non-worsening, by taking a weighted convex combination of the gradients of the individual objectives $\nabla f_0$ and $\nabla f_1$,
\begin{equation}
    d(\mathbf{x}) = w_{0}\cdot\nabla f_0(\mathbf{x}) + w_{1}\cdot\nabla f_{1}(\mathbf{x}),
    \label{eq:weightedsum}
\end{equation}
with $w_{0}, w_{1} \geq 0$. 
In general, there exist infinitely many search directions $d(\mathbf{x})$ for which all objectives are non-worsening, and different methods have been developed in which a single descent direction is computed \cite{fliege2000,schaffler02,desideri2012multiple}.
While this provides an approach to converge to Pareto optimal solutions, it does not tell us directly how to take solution diversity into account, which has shown to be non-trivial \cite{bosman2011gradients,schutze2016directed}. 
We therefore consider a different avenue to handle gradients for MO optimization in this work, which is to cast the MO problem as a single-objective (SO) optimization problem, in which a quality indicator is used to quantify the quality of an approximation set \cite{Emmerich2018MOtutorial,li2019quality}. One popular quality indicator is the hypervolume indicator \cite{zitzler1998multiobjective}, which measures the volume in objective space that is dominated by an approximation set. The hypervolume indicator is currently the only known indicator that is Pareto-compliant, meaning that solutions in a set with maximal hypervolume are Pareto optimal \cite{fleischer2003}, and it furthermore takes diversity intrinsically into account \cite{auger2009HV}. Additionally, the hypervolume indicator is differentiable with respect to a problem's objective functions in strictly non-dominated points
which allows determining gradient weights via the chain rule \cite{emmerich2014time}.

A limitation of the hypervolume indicator however is that it ignores dominated solutions. This prevents the use of the hypervolume indicator directly in indicator-based MO optimization, as it cannot be used to steer dominated solutions to a non-dominated area in the search space \cite{toure2019uncrowded}. SMS-EMOA \cite{beume2007sms} overcomes this limitation by using non-dominated sorting to create subsets of solutions such that solutions within a subset are non-dominated. Consequently, each solution's hypervolume contribution with respect to its subset can be computed and used to steer the solution towards the Pareto front. The hypervolume indicator gradient ascent multi-objective optimization (HIGA-MO) algorithm \cite{wang2017hypervolume} computes hypervolume gradients for solutions in subsets created by non-dominated sorting and thus achieves gradient-based steering for dominated solutions.
An approach to incorporate dominated solutions into a hypervolume-based indicator is the uncrowded hypervolume improvement \cite{toure2019uncrowded} which was combined with the newly presented Sofomore framework to perform optimization by interleaving single-objective optimizers. In \cite{maree2020}, this quality indicator for single solutions was recently converted into a quality measure for solution sets, called the uncrowded hypervolume (UHV), which is directly suitable for indicator-based MO optimization. The resulting UHV problem was then efficiently solved with the gene-pool optimal mixing evolutionary algorithm by exploiting UHV-specific properties (UHV-GOMEA). Earlier, \cite{emmerich2007gradient} had already described directing dominated solutions to the approximation boundary by minimizing their Euclidean distance to the boundary in the context of gradient ascent for HV optimization. 

In this work, we formulate gradient expressions for the UHV, such that it can be used directly in SO gradient ascent schemes. (Note that the UHV needs to be maximized, independent of whether the underlying MO problem is a minimization or maximization problem.) To demonstrate this, we solve it with the same scheme as used by HIGA-MO, and with \textit{Adam}, a preeminent method for efficient stochastic optimization \cite{kingma2014adam}. We further compare UHV gradient ascent to HIGA-MO, and the EAs UHV-GOMEA and Sofomore-GOMEA \cite{maree2020}. For the experimental comparison, we employ simple quadratic benchmark functions similar to benchmarks used in \cite{maree2020} and also the Walking Fish Group (WFG) benchmark set \cite{huband2005scalable}. Additionally, for a fair comparison to EAs, we study the performance of the gradient-based methods in a black-box setting, where gradient information of the MO problem is not available, by using a finite difference gradient approximation.
The remainder of this paper is organized as follows. In Section~\ref{sec:prelim}, we introduce preliminaries of the (uncrowded) hypervolume indicator. In Section~\ref{sec:UHV-Grad}, we introduce our UHV gradient ascent algorithm. Experimental comparisons are described in Section~\ref{sec:experiments}, followed by a discussion in Section~\ref{sec:discussion}.


\section{Uncrowded hypervolume optimization}
\label{sec:prelim}
We consider MO problems given by a to-be-minimized $m$-dimensional objective function $\mathbf{f} : \mathcal{X} \rightarrow \mathbb{R}^m$, where $\mathcal{X}\subseteq\mathbb{R}^n$ is the $n$-dimensional (box) constrained decision space. We focus on the bi-objective case $m = 2$ in this work. Let $\mathbf{x}\in\mathcal{X}\subseteq\mathbb{R}^n$ be a solution of the MO problem, which we from now on refer to as an \textit{MO-solution}. The goal of MO optimization is to obtain a set of (near-)Pareto-optimal MO-solutions $\calS\subset\mathcal{X}$ of manageable size. To evaluate the quality of $\calS$, we use the uncrowded hypervolume (UHV) indicator \cite{maree2020}, which measures the area in objective space enclosed by the non-dominated MO-solutions in $\calS$ and a reference point $r = (r_0, r_1)$ (as the hypervolume indicator \cite{zitzler1998multiobjective}), and uses the uncrowded distance \cite{toure2019uncrowded} (explained below) to steer dominated MO-solutions. As $\mathcal{S}$ can contain dominated MO-solutions, let $\mathcal{A}$ be the approximation set of $\mathcal{S}$, i.e., the largest subset of $\mathcal{S}$ that contains only non-dominated MO-solutions.

\begin{figure}[t]
     \centering
    \includegraphics[width = \textwidth]{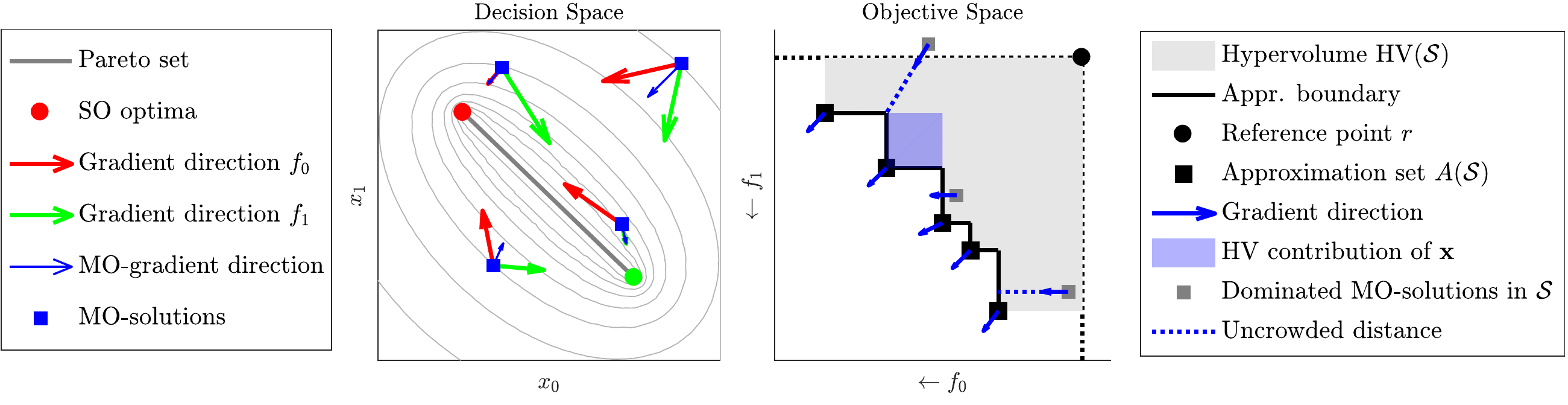}
    \vspace*{-0.8cm}
     \caption{Illustration of UHV gradient ascent on a bi-objective problem. The MO-gradient direction in decision space (left subfigure) is a weighted linear combination of the SO gradients, where the weights are determined based on the UHV gradient direction in objective space (right subfigure).}
     \label{fig:uhv_grad}
     \vspace*{-0.5cm}
\end{figure}

In order to search the space of solution sets, $\raisebox{.15\baselineskip}{\Large\ensuremath{\wp}}(\mathcal{X})$, a parameterization of solution sets is required \cite{wang2017hypervolume,Beume07,maree2020}. For this, we consider sets $\calS_p$ of a fixed size of $p$ MO-solutions, and simply concatenate the decision variables of all MO-solutions into a single vector $X\in\mathbb{R}^{np}$, i.e., $X = [ \mathbf{x}_0 \; \cdots \; \mathbf{x}_{p-1} ]$, similar to notation used in \cite{emmerich2014time}. Additionally, let $Y\in\mathbb{R}^{p\times m}$ be the matrix of concatenated objectives values corresponding to $X$, i.e., $Y_{i,0:m-1} = \mathbf{y}_i = \mathbf{f}(\mathbf{x}_i)$. Finally, let $F : \mathbb{R}^{np} \rightarrow \mathbb{R}^{m\times p}$ be the operator that evaluates the entire solution set given by $X$, i.e., $Y = F(X)$. This implies that an evaluation of $F$ consists of $p$ evaluations of the MO problem (MO-evaluations). The resulting SO UHV-based optimization problem can then be formulated as,
\begin{equation}
\label{eqn:ppsn20_UHVmop}
\begin{split}
\text{maximize} \quad & g(X) = \text{UHV}(F(X)) 
= \text{HV}(F(X)) - \mbox{UD}(F(X)), 
\\
\text{with} \quad & \mathbf{f} : \mathcal{X} \subseteq \mathbb{R}^n \rightarrow \mathbb{R}^m, \quad F : \mathbb{R}^{np} \rightarrow \mathbb{R}^{m\times p}, \quad X \in \mathbb{R}^{np},
\end{split}
\end{equation}
where $\mbox{HV} : \mathbb{R}^{m\times p} \rightarrow \mathbb{R}_{\geq 0}$ is the hypervolume indicator \cite{zitzler1998multiobjective} and $\mbox{UD} : \mathbb{R}^{m\times p} \rightarrow \mathbb{R}_{\geq 0}$ is the mean of the uncrowded distances $\mbox{ud}(\mathbf{y},Y)$ \cite{toure2019uncrowded}, which measure the shortest distance of a point $\mathbf{y}$ towards the domination boundary of $Y$ in objective space. It is called the uncrowded distance as the nearest point on the boundary of $Y$ is generally away from points in $Y$. The UD is then given by,
\begin{equation}
\label{eqn:ud}
    \mbox{UD}(Y) = \frac{1}{p} \sum_{i = 0}^{p-1} \text{ud}(\mathbf{y}_i,Y)^m.
\end{equation}
We refrain from repeating a mathematical definition here, but provide an illustration in Figure~\ref{fig:uhv_grad}. Note that, in contrast to \cite{maree2020}, we only consider the interior boundary of $Y$ here, which was found to improve performance in preliminary experiments, as the extreme points of $Y$ are often already well-positioned close to the extremes of the approximation front (i.e. the approximation set in objective space), and steering additional points into the same location causes undesired computational overhead. Finally, note that the UHV is equivalent to the hypervolume indicator when all MO-solutions in $\mathcal{S}$ are non-dominated, which implies that the UHV is still Pareto-compliant on the space of approximation sets.

\section{UHV gradient ascent}
\label{sec:UHV-Grad}
We apply a gradient ascent scheme to $g(X) = \mbox{UHV}(F(X))$ in Eqn.~\eqref{eqn:ppsn20_UHVmop}. For this, we use the gradient of the hypervolume indicator as was derived in \cite{emmerich2014time}. We briefly describe the concept here, but refer the reader to \cite{emmerich2014time} for a rigorous mathematical derivation and analysis. 
The gradient $\nabla g(X) = \nabla \mbox{UHV}(F(X))$ can be split up into subvectors corresponding to different MO-solutions by using that $X = [ \mathbf{x}_0 \; \cdots \; \mathbf{x}_{p-1} ]\in\mathbb{R}^{np}$,
\begin{equation}
    \nabla g(X) = \dd{\mbox{UHV}(F(X))}{X} = \left[ \dd{\mbox{UHV}(F(X))}{\mathbf{x}_0} \; \cdots \; \dd{\mbox{UHV}(F(X))}{\mathbf{x}_{p-1}} \right].
\end{equation}
We now apply the chain rule to each of the subvectors $i$ by using $\mathbf{y}_i = \mathbf{f}(\mathbf{x}_i)$,
\begin{equation}
    \dd{\mbox{UHV}(F(X))}{\mathbf{x}_i} = \dd{\mbox{UHV}(F(X))}{F(X)} \cdot \dd{F(X)}{\mathbf{x}_i} = \sum_{j = 0}^{p-1} \dd{\mbox{UHV}(F(X))}{\mathbf{y}_j} \cdot \dd{\mathbf{y}_j}{\mathbf{x}_i},
\end{equation}
where we can now use that $\dd{\mathbf{y}_j}{\mathbf{x}_i} = \mathbf{0}$ for $j \neq i$, as the fitness values of $\mathbf{y}_j = \mathbf{f}(\mathbf{x}_j)$ do not depend on $\mathbf{x}_i$. For $j = i$, we have $\dd{\mathbf{y}_i}{\mathbf{x}_i} = [ \nabla f_0(\mathbf{x}_i) \; \nabla f_1(\mathbf{x}_i)]$, which are simply the gradients of the MO problem. This gives,
\begin{equation}
    \dd{\mbox{UHV}(F(X))}{\mathbf{x}_i} = \dd{\text{UHV}(F(X))}{f_0(\mathbf{x}_i)} \cdot \nabla f_0(\mathbf{x}_i) + \dd{\mbox{UHV}(F(X))}{f_1(\mathbf{x}_i)} \cdot \nabla f_1(\mathbf{x}_i).
\end{equation}
Note the correspondence of this expression with the weighted search direction in Eqn.~\eqref{eq:weightedsum}. 
Directly using the objective space gradients to determine the search direction would cause MO-solutions that contribute more to the UHV to make big steps, and MO-solutions that contribute little to slowly \textit{creep}, which was noted earlier \cite{hernandez2014hypervolume,wang2017hypervolume}. To overcome this, we normalize the objective gradients by setting $W = \left\lVert \left[ \dd{\text{UHV}}{f_0(\mathbf{x}_i)} \; \dd{\text{UHV}}{f_1(\mathbf{x}_i)} \right]\right\rVert$, which gives us the desired search direction,
\begin{equation}
    \frac{1}{W}\dd{\mbox{UHV}(F(X))}{\mathbf{x}_i}\! =\! \frac{1}{W}\dd{\text{UHV}(F(X))}{f_0(\mathbf{x}_i)} \cdot\! \nabla f_0(\mathbf{x}_i) + \frac{1}{W}\dd{\mbox{UHV}(F(X))}{f_1(\mathbf{x}_i)} \cdot \nabla f_1(\mathbf{x}_i).
\end{equation}
It now remains to find an expression for the objective space gradients. We now use that $\mbox{UHV} = \mbox{HV} - \mbox{UD}$. For both objectives $k = \{0,1\}$, this gives, $$\dd{\text{UHV}(F(X))}{f_k(\mathbf{x}_i)} = \dd{\text{HV}(F(X))}{f_k(\mathbf{x}_i)} - \dd{\text{UD}(F(X))}{f_k(\mathbf{x}_i)}.$$
Whenever $\mathbf{x}_i$ is a dominated MO-solution, it has no contribution to the hypervolume, and the first term is therefore equal to zero. For the second term, let $\mathbf{s}(\mathbf{f}(\mathbf{x}_i))\in\mathbb{R}^m$ be the point towards which the uncrowded distance is computed, i.e., the nearest point to $\mathbf{f}(\mathbf{x}_i)$ on the approximation boundary given, as illustrated in Figure~\ref{fig:uhv_grad}. Using the definition of UD in Eqn.~\eqref{eqn:ud}, we obtain the final expression for objective-space derivative for dominated MO-solutions, $$\dd{\text{UD}(F(X))}{f_k(\mathbf{x}_i)} =  \frac{1}{p} \dd{}{f_k(\mathbf{x}_i)} \lVert \mathbf{f}(\mathbf{x}_i) - \mathbf{s}(\mathbf{f}(\mathbf{x}_i)) \rVert^m .$$

Whenever $\mathbf{x}_i$ is a non-dominated MO-solution, the objective-space hypervolume gradient can be computed by the approach described in \cite{emmerich2014time}. Conceptually, the computation can be reduced to the objective-space gradient of the hypervolume contribution of that MO-solution, which is easily computed when the neighbouring MO-solutions on the approximation front are known (Figure~\ref{fig:uhv_grad}). Additionally, whenever $\mathbf{x}_i$ is a non-dominated MO-solution, it determines the approximation boundary, which is used in the computation of the UD for other MO-solutions. Therefore, $\dd{\text{UD}(F(X))}{f_k(\mathbf{x}_i)}$ is potentially non-zero. In that case, the UD can be improved at the cost of worsening non-dominated MO-solutions, as this reduces the uncrowded distance of dominated MO-solutions. This is undesirable, and we therefore explicitly set $\dd{\text{UD}(F(X))}{f_k(\mathbf{x}_i)} = 0$ for non-dominated $\mathbf{x}_i$, although preliminary experiments showed that performance is largely unaffected by this. Finally, we consider the case in which $\mathbf{x}_i$ is \textit{weakly dominated}, which occurs for pairs of MO-solutions with at least one coinciding objective value. In this case, the objective-space gradient of the UHV is undefined \cite[Proposition 3]{emmerich2014time}. To prevent such case, we consider these points to be strongly dominated, and (temporarily) worsen the objective value(s) that are shared with other MO-solutions by a small value $\varepsilon$, which allows us to compute the uncrowded distance as before. Since objective space gradients are normalized, the actual choice of $\varepsilon$ is irrelevant as long as it is small enough so that the weakly dominated MO-solution does not get dominated by other MO-solutions.

\begin{table}[t]
\caption{UHV gradient ascent schemes for maximizing $g(X)$.}
\vspace*{-0.3cm}
\label{tab:GA_schemes}

    \begin{tabular}{p{0.48\textwidth}|p{0.5\textwidth}}
    \toprule
         \textbf{Adam} \cite{kingma2014adam} &  \textbf{GA-MO} \cite{higamogithub} \\
         \midrule
         \textbf{Initial values}: $\gamma^0 = \lVert\mathcal{X}_\text{init}\rVert\cdot10^{-2}$, $b_0 = 0.9$, $b_1 = 0.999$, $b_2 = 0.99$, $\varepsilon = 10^{-16}$, $\mathbf{m}^{-1} = \mathbf{v}^{-1} = \mathbf{0}$.
         & 
         \textbf{Initial values}:  $c = 0.1$, $\alpha = 0.7$, $\beta = 0.7$, and for $i = 0,\ldots,(p-1)$: $\gamma^{-1}_i = \lVert\mathcal{X}_\text{init}\rVert\cdot10^{-2}$, $\mathbf{n}^{-1}_i = \mathbf{0}$, $m_i^{-1} = 0$. \\
         \midrule

For $t = 0,1,\ldots,$
{$$\setlength{\jot}{0pt} 
\!\begin{aligned}
        \mathbf{m}^t & = b_0 \mathbf{m}^{t-1} + (1-b_0)\nabla g(X^t), \\
        \mathbf{v}^t & = b_1 \mathbf{v}^{t-1} + (1-b_1)\nabla g^2(X^t), \\
        X^{t+1} & = X^t + \gamma^t \frac{\mathbf{m}^t / (1-(b_0)^{t+1})}{\sqrt{\mathbf{v}^t / (1-(b_1)^{t+1})} + \varepsilon},\\
        \gamma^{t+1} & = \begin{cases} b_2 \gamma^t, & \mbox{if } g(X^{t+1}) \leq g(X^t), \\
        \gamma^t, & \mbox{else. }
        \end{cases}
    \end{aligned}$$}
         
          & 

For $t = 0,1,\ldots,$
{$$\setlength{\jot}{0pt} 
\scriptsize
\!\begin{aligned}
\mbox{d}^-  & = \min_{l,k\in\{0,\ldots,(p-1)\}, l\neq k}\lVert\mathbf{x}^t_l - \mathbf{x}^t_k\rVert, \\
\mbox{d}^+ & = \max_{l,k\in\{0,\ldots,(p-1)\}, l\neq k}\lVert\mathbf{x}^t_l - \mathbf{x}^t_k\rVert, \\
\gamma^\text{UB} & = \beta (\mbox{d}^+ + \mbox{d}^-) / 2. \\ 
\end{aligned}
$$}
\; \; \; For $i = 0,\ldots,(p-1)$,
{$$\setlength{\jot}{0pt} 
\!\begin{aligned}
        \mathbf{n}_i^t & = \nabla g(\mathbf{x}_i^t) /  \lVert \nabla g(\mathbf{x}_i^t) \rVert,  \\
        m^t_i & = (1-c)m_i^{t-1} + c \langle \mathbf{n}_i^{t-1},\mathbf{n}_i^{t}\rangle, \\
        \gamma_i^t & = \min\{\gamma^\text{UB}, \gamma_i^{t-1} e^{\alpha m_i^t}\}, \\
        \; \; \; \; \; \; \; \; \; \mathbf{x}_i^{t+1} & = \mathbf{x}_i^t + \gamma_i^t \mathbf{n}_i^t. 
    \end{aligned}
$$} \vspace*{-0.4cm}
    \\ \bottomrule
    \end{tabular}
    \vspace*{-0.4cm}
\end{table}

\subsection{Gradient ascent schemes}
We use two gradient ascent schemes for UHV gradient ascent, as listed in Table~\ref{tab:GA_schemes}. The first scheme we consider is \textit{Adam} \cite{kingma2014adam} (UHV-Adam), which is a popular method for stochastic gradient descent. Adam uses a variance-corrected weighted average of current and previous gradients. In contrast to the original formulation, we set $\varepsilon$ to machine precision, and we add a very simple step size shrinking scheme in which the step size is reduced if no improvement was found. The second scheme is the GA-MO scheme (UHV-GA-MO) used in the Python implementation of HIGA-MO \cite{higamogithub}. GA-MO updates the step size for each MO-solution separately using a weighted average of search directions' inner products as input for an exponential cooling scheme. We adapted the weight used for averaging inner products $c = 0.1$ (from $c = 0.2$) as both HIGA-MO and UHV-GA-MO showed stagnation in preliminary experiments with $c = 0.2$. Additionally, we changed the upper bound on the step size $\gamma^\text{UB}$ to be also based on $d^+$, the maximum distance between two MO-solutions in decision space, as for the UHV objective function, two dominated MO-solutions could be steered to the same point on the front, and only basing it on the minimum distance $d^-$ could shrink $\gamma$ prematurely. For both schemes, we use projected gradients (i.e., boundary repair) to handle box-constrained search spaces. Initial MO-solutions are initialized uniformly random in a box $\mathcal{X}_\text{init} \subseteq \mathcal{X}$, and the initial step size is based on the maximum initialization range in any dimension, which we denote by $\lVert\mathcal{X}_\text{init}\rVert$. Implementations of UHV-Adam and UHV-GA-MO are available at \url{https://github.com/scmaree/uncrowded-hypervolume}.

\subsection{Finite difference gradient approximation}
To assess the performance of gradient-based algorithms in a black-box scenario, where exact gradients are not known, finite forward difference gradient approximations (FD) are used. The FD step size is set to $h = 10^{-6}\cdot \bar{\gamma}_t$, where $\bar{\gamma} = \sum_{i=0}^{p-1} \gamma_i^t$ for UHV-GA-MO, and $\bar{\gamma} = \gamma^{t}$ for UHV-Adam.
In this way, $h$ is always smaller than the mean step size. If the FD step violates the search space's box-constraints, backward differences are used. 
Estimating both objectives' gradients in one MO-solution requires $n$ additional MO-evaluations, the number of MO-evaluations thus increases from $p$ to $(1+n)\cdot p$ per iteration. When using FD, we refer to our methods as UHV-Adam-FD and UHV-GA-MO-FD.


\section{Experiments}
\label{sec:experiments}
Experiments are conducted on several bi-objective problems: four bi-objective problems with known gradients as defined in Table~\ref{tab:benchmarks} and nine box-constrained problems from the WFG benchmark suite \cite{huband2005scalable,wfgcpp,wessing2018optproblems}. For each algorithm, the \textit{best} approximation set obtained so far is recorded over the run of that algorithm, where quality is measured by the algorithm itself, i.e., based on the HV or UHV. Performance is measured by the number of MO function evaluations (MO-evaluations), where we define one MO-evaluation as the computation of $f_0, f_1, \nabla f_0, $ and $\nabla f_1$ at once. Note that the evaluation of $X$, which models a solution set $\calS_p$ of size $p$, therefore costs $p$ MO-evaluations.  All problems are run with a fixed hypervolume reference point $r = (11,11)$, which is rather far away from the Pareto front, as this puts additional importance towards obtaining the end points of the front \cite{auger2009HV}. However, even with this choice of reference point, the endpoints are not always included in the approximation set with optimal hypervolume, depending on the shape of the front \cite{auger2009HV}. 

We compare the two UHV gradient ascent schemes, UHV-Adam and UHV-GA-MO, to the EAs UHV-GOMEA-Lt and Sofomore-GOMEA from \cite{maree2020}. We furthermore consider the gradient-based HIGA-MO. UHV-GOMEA-Lt uses a linkage tree in which at most a few MO-solutions are updated simultaneously. A full description of UHV-GOMEA-Lt and Sofomore-GOMEA can be found here \cite{maree2020}. 
We used the Python implementation of HIGA-MO \cite{higamogithub}, but extended it with a dynamic reference point so that hypervolume gradients can also be computed for solution sets with $f_0^+ > r_0$ or $f_1^+>r_1$ (which is not an issue for the UHV-based algorithms), where $f_i^+$ is the worst value for the $i^\text{th}$ objective in the set. The dynamic reference point is set to $\hat{r} = (\max\{1.1f_0^+, r_0\},$ $\max\{1.1f_1^+, r_1\})$. Algorithmic performance is always evaluated with respect to $r$.

As performance indicators, we consider the difference with the optimal hypervolume (for $p$ MO-solutions) $\Delta\mbox{HV}(\mathcal{A}_p) = \mbox{HV}(\mathcal{A}^\star_p) - \mbox{HV}(\mathcal{A}_p)$, where $\mathcal{A}_p = A(\calS_p)$ is the approximation set given by $\calS_p$, and $\mathcal{A}^\star_p$ is the approximation set with optimal hypervolume. The second measure we consider is the generational distance (GD) \cite{Zitzler2003}.  
The GD for problems 0 and 2 is computed analytically from their known Pareto set \cite{maree2020}.
For problems 1 and 3, the GD is computed based on a sample of 5000 MO-solutions from a reference set. The GD is not Pareto compliant, but it is a useful tool to measure proximity to the Pareto set. Finally, we consider $|\mathcal{A}_p|$, i.e., the number of non-dominated MO-solutions in $\calS_p$, which we use to measure how well different mechanisms for handling dominated MO-solutions perform. Unless mentioned otherwise, all experiments are repeated 10 times, and medians and inter-quartile ranges (IQR) are shown.

\begin{table}[t]
\caption{Quadratic bi-objective benchmark problems. $R$ applies a $45^{\circ}$ rotation along all axes. All MO-solutions are initialized in $[-2,2]^{n \times p}$ for problems 0--2 and $[0,2]^{n \times p}$ for Problem 3.
}
\label{tab:benchmarks}
\vspace*{-0.2cm}
\scriptsize
\begin{adjustbox}{width=\textwidth,center}
\begin{tabular}{cccccc}
    \# & Problem & $f_{0}$  & $f_{1}$ & Properties\\
    \midrule
    0
    &
    \begin{tabular}{c}
    convex \\ bi-sphere
    \end{tabular}&
    $f(\textbf{x}) = \sum_{i=0}^{n-1}(x_{i})^{2}$  & 
    \begin{tabular}{c}
    $f(\textbf{x} - \textbf{c}),$ with $\textbf{c} = [1 \; 0 \dots 0]$
    \end{tabular}
    & 
    \begin{tabular}{c}
    decomposable,
    \end{tabular}
    \\
    \midrule
    1
    &
    \begin{tabular}{c}
    sphere \&\\ rotated ellipsoid
    \end{tabular} & 
    $\frac{1}{n}f(\mathbf{x})$& 
    \begin{tabular}{c}
    $(\sqrt{W}R\mathbf{x} - \sqrt{W}c)^{\intercal}(\sqrt{W}R\mathbf{x} - \sqrt{W}\mathbf{c})$, 
    $W_{i,i} = 10^{\frac{-6i}{n-1}}$
    \end{tabular}
    & 
    \begin{tabular}{c}
    non-decomposable,
    \\
    ill-conditioned
    \end{tabular}
    \\
    \midrule
    2
    &
    \begin{tabular}{c}
    concave \\ bi-sphere
    \end{tabular} & 
    $f(\mathbf{x})^{\frac{1}{4}}$& 
    $f(\mathbf{x}-\mathbf{c})^{\frac{1}{4}}$&
    \begin{tabular}{c}
    decomposable,
    \\
    concave front
    \end{tabular}
    \\
    \midrule
    3
    &
    \begin{tabular}{c}
    sphere \&\\ Rosenbrock
    \end{tabular} & 
    $f(\mathbf{x})$ &
    $\frac{1}{(n-1)} \sum_{i=0}^{(n-2)} 100\left(x_{(i+1)}- x_{i}^{2}\right)^{2}+\left(1-x_{i}\right)^{2}$ & 
    \begin{tabular}{c}
    bimodal,
    \\
    chained dependencies
    \end{tabular}
     \\

\end{tabular}
\end{adjustbox}
\vspace*{-0.4cm}
\end{table}

\subsection{Convergence in hypervolume on the quadratic functions}
\label{sec:exp1}
\begin{figure}[t]
    \vspace*{-0.0cm}
    \centering
    \includegraphics[width = \textwidth]{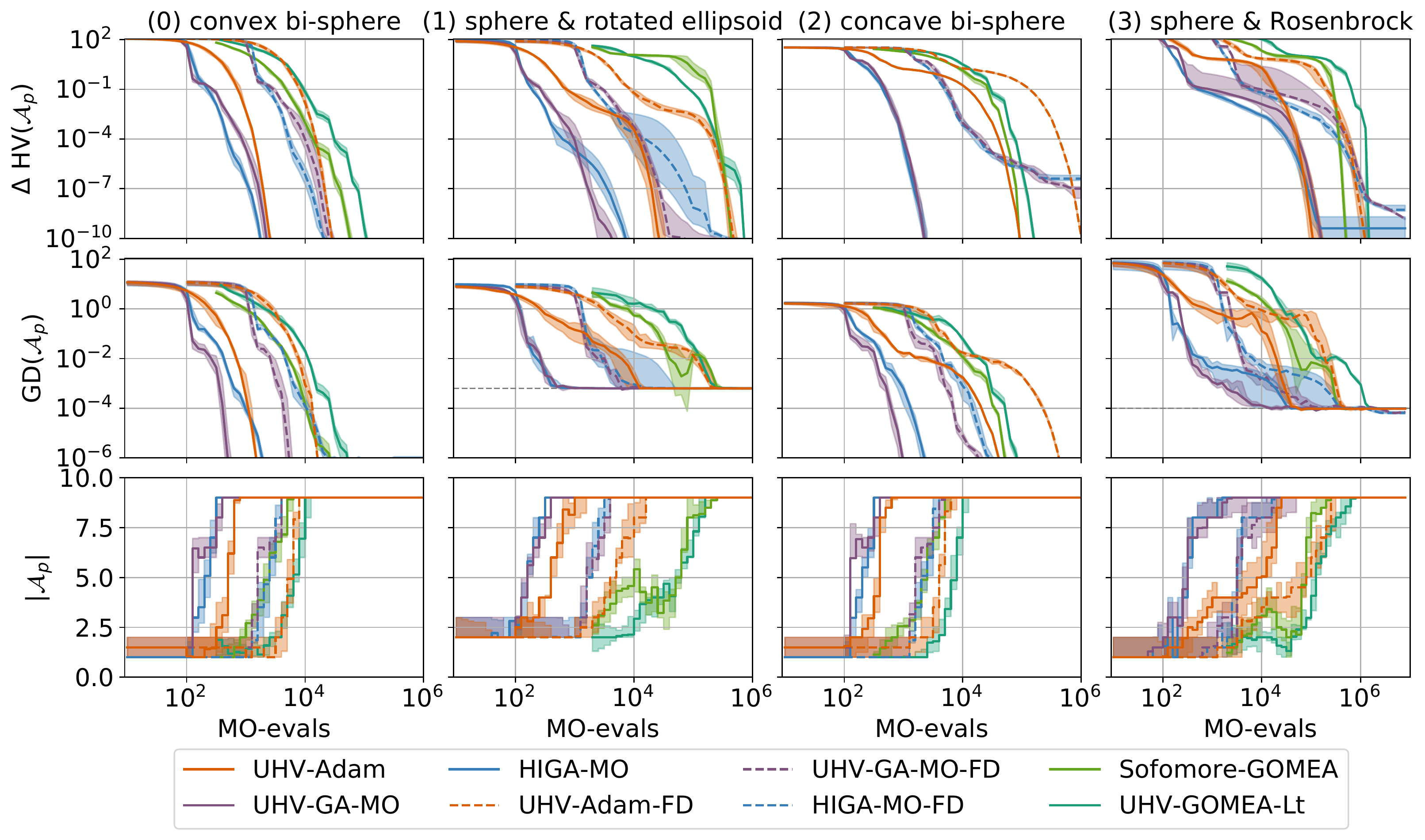}
    \vspace*{-0.8cm}
    \caption{Results for the different algorithms with $p = 9$ on the benchmark problems in Table~\ref{tab:benchmarks} with $n = 10$. Lines indicate median values and shaded areas represent the IQR. Dashed lines correspond to gradient-based algorithms with finite difference gradient approximations. MO-evals: MO-evaluations.}
    \label{fig:exp1}
    \vspace*{-0.7cm}
\end{figure}

For the first experiment, we consider the quadratic bi-objective functions from Table~\ref{tab:benchmarks}, with problem dimensionality $n = 10$. We consider solution sets $\calS_p$ of size $p = 9$. As we do not know $\mbox{HV}(\mathcal{A}^\star_p)$ analytically, the target HV is set to maximal HVs obtained from lower dimensional instances.
For UHV-Adam, UHV-GA-MO, and HIGA-MO, the initial step size $\gamma^{0}$ is set to one percent of the mean initialization range. As in \cite{maree2020}, UHV-GOMEA-Lt and Sofomore-GOMEA are run with population size $N = 31$ for the decomposable problems 0 and 2, and $N=200$ otherwise. All optimizers are run for $10^{6}$ MO-evaluations or until convergence criteria are met.
Results are shown in Figure~\ref{fig:exp1}. All gradient-based algorithms reach the target hypervolume in all problems except for HIGA-MO on Problem 3 which converges close to the target HV. HIGA-MO's performance is more volatile across runs which is especially visible in problems 1 and 3. Both EAs, UHV-GOMEA-Lt and Sofomore-GOMEA, always obtain the target hypervolume, but require substantially more MO-evaluations. The IQRs of gradient-based and EA-based algorithms only rarely intersect, indicating that the faster convergence of gradient-based algorithms is robust to random initialization. These differences in performance are also reflected in $\GDA$. All gradient-based algorithms obtain $|\mathcal{A}_p|= 9$ non-dominated MO-solutions faster than EAs. UHV-GA-MO and HIGA-MO reach $|\mathcal{A}_p|= 9$ sooner than UHV-Adam which indicates that the gradient ascent scheme (i.e., GA-MO vs. Adam) has a larger effect on quickly finding non-dominated MO-solutions than the strategy for handling dominated MO-solutions (i.e., UD vs. non-dominated sorting of MO-solutions into multiple fronts in HIGA-MO). 

\subsubsection{Finite difference gradient approximation.}
In practice it may well happen that analytic gradients are not available.
In Figure~\ref{fig:exp1}, it can be seen that the gradient-based algorithms lose much of their advantage over EA-based algorithms when relying on finite difference gradient approximations. Their convergence is only slightly faster (Problem 0), similar (Problem 3), or EAs now clearly outperform them (Problem 2). In Problem 2, HIGA-MO-FD converges prematurely. Only in Problem 1, there is still evidence of advantages for gradient-based algorithms: UHV-GA-MO-FD still convergences more than 10 times faster than EA-based algorithms. This problem is highly dependent and ill-conditioned, and a large population is required for the EAs to solve this problem, while gradient-based algorithms directly capture these dependencies. These results also show that finite difference gradient approximations not only increase the computational cost per iteration, but could also worsen convergence rates or even cause stagnation. This is especially true for the GA-MO scheme. UHV-Adam-FD does however not show a deterioration in the rate of convergence (besides the expected shift of a factor $1+n$). Adam was developed for stochastic gradient descent, and uses a weighted average of current and past gradients instead of the gradient itself, enabling it to handle the imprecise gradient approximations.

\subsection{Effect of the number of MO-solutions $p$}
When $p$ is increased, the size of the to-be-optimized approximation set $\mathcal{A}_p$ is also increased. This makes the resulting UHV optimization problem more difficult, and dependency modelling becomes essential in order to obtain the optimal distribution of MO-solutions along the front with UHV-GOMEA-Lt \cite{maree2020}. 
To investigate the dependence of convergence speed on the number of MO-solutions $p$, all gradient-based optimizers are applied on problems 0--2 with $p = 2^j + 1$ for $j = 1,\ldots,7$  and $n=10$. Problem 3 is excluded as premature convergence to its local optimum would obfuscate the comparison.
All optimizers are run for $10^{7}$ MO-evaluations or until convergence criteria are met. The target HV is set to the maximal HV found across all algorithms.
Parameter settings ($N$ for UHV-GOMEA-Lt, $\gamma^{0}$ otherwise) were tuned experimentally across all $p = 2^j + 1$: $\gamma^{0}$ is set to $4\cdot 10^{-2}$ for UHV-Adam, $4\cdot 10^{-4}$ for UHV-GA-MO, and to $4\cdot 10^{-3}$ for HIGA-MO.  
UHV-GOMEA-Lt's population sizes are scaled in $p$ as the larger parameter spaces require larger populations: $N(p) = \lceil 0.76 p^{\frac{1}{4}}N_{\text{base}}\rfloor$, where $\lceil \cdot \rfloor$ is the rounding operator. $N_{\text{base}}=31$ for problems 0 and 2 and $N_{\text{base}} = 200$ for Problem 1 as in Section~\ref{sec:exp1}.
Sofomore-GOMEA interleaves optimizations of individual MO-solutions, therefore $N$ does not need to be scaled in $p$ and $N$ is set to $N_{\text{base}}$.
\begin{figure}[t]
    \centering
    \includegraphics[width = \textwidth]{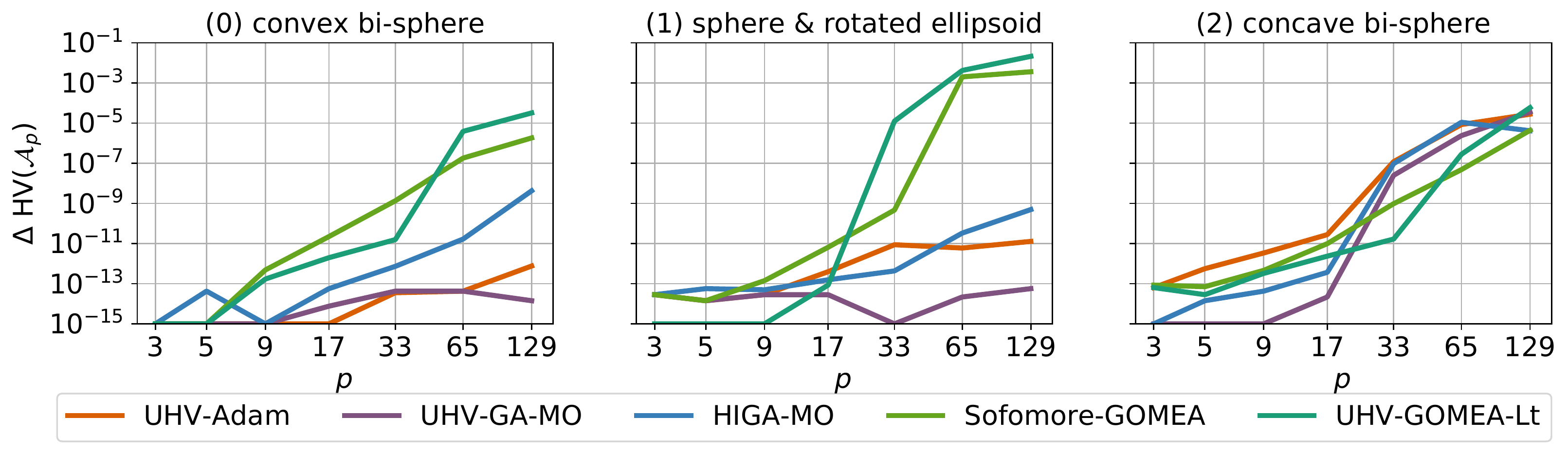}
    \caption{The median distance to the target $\mathrm{HV}$ after $10^{7}$ MO-evaluations of all algorithms for problems 0--2 with varying $p$ and $n=10$ over 10 repetitions.}
    \label{fig:exp_maxhv_best}
    \vspace*{-0.6cm}
\end{figure}
The median $\DHV$ in problems 0--2 with varying $p$ is shown in Figure~\ref{fig:exp_maxhv_best}. All algorithms always reach the target HV with $10^{-10}$ accuracy for $p\leq 17$. As $p$ increases, the $\DHV$ of UHV-GOMEA-Lt and Sofomore-GOMEA increases across problems. All gradient-based algorithms obtain lower $\DHV$ values than both EAs as $p$ increases with the exception of Problem 2, in which UHV-GOMEA-Lt and Sofomore-GOMEA scale better in $p$.

\subsection{WFG benchmark}
The WFG test suite \cite{huband2005scalable} consists of 9 benchmark functions with different properties. WFG1 is decomposable, but has a flat region in the decision space, which could cause stagnation. WFG2, WFG4, and WFG9 have one or more multimodal objectives, which are expected to be difficult for gradient-based algorithms. Problems WFG4--9 have concave fronts, WFG1 has a convex front, WFG2 has a disconnected convex front, and WFG3 has a linear front. We use finite difference approximations for the gradient-based algorithms.
We again consider bi-objective problems, and use $k_\text{WFG} = 4$ position variables and $l_\text{WFG} = 20$ distance variables, resulting in a total of $n = 24$ decision variables as originally chosen in \cite{huband2006review}. We solve these benchmarks with approximation sets of size $p = 9$ and a limited computational budget of $10^5$ MO-evaluations. All experiments are repeated 30 times. Differences are tested for statistical significance (up to 4 decimals) by a Wilcoxon rank sum test with $\alpha = 0.05$, pairwise to the best. Ranks (in brackets) are computed based on the mean hypervolume. All statistics are computed per table.  For the gradient-based algorithms, we set $\gamma^{0} = \lVert\mathcal{X}_\text{init}\rVert \cdot 10^{-2}$, and a population size of $N = 200$ was used for the population-based algorithms.
Results on the WFG benchmark are shown in Table~\ref{tab:ppsn20_grad_wfg}. UHV-Adam-FD performs best overall, while UHV-GA-MO-FD has worse performance on most problems. As expected, the gradient-based algorithms perform worse on the multi-modal problems WFG4 and WFG9. WFG2 has only one multimodal objective which does not seem to be a problem for UHV-Adam-FD. All algorithms have difficulties with the flat region in WFG1, and the worst overall hypervolume values are obtained for this problem. The only algorithm that has an explicit mechanism for handling flatness is HIGA-MO-FD, which consequently performs best for WFG1. HIGA-MO-FD re-initializes MO-solutions if the gradient is zero, which does not help traversing the plateau, but increases diversity. Note that we did not add any mechanism to handle flatness in the other algorithms. Especially with gradient-based algorithms, flatness is easily detected and a mechanism could be added to improve performance. On the WFG problems with concave fronts, UHV-Adam-FD performs well, in contrast to the previous results on the concave bi-sphere, showing that it does not have difficulties with concavity in general.

\begin{table}[htp]
\caption{Results on the WFG Benchmark. Hypervolume values are shown (mean, $\pm$ standard deviation (rank)). Finite differences (FD) are used for the gradient-based algorithms. Bold are best scores per problem, or those not statistically different from it.}
\vspace*{-0.3cm}
\label{tab:ppsn20_grad_wfg}
\begin{adjustbox}{width=\textwidth,center}


\begin{tabular}{crrrrrrr}
\toprule
Problem & Sofomore-GOMEA & UHV-GOMEA-Lt & UHV-ADAM-FD & UHV-GA-MO-FD & HIGA-MO-FD \\
 \midrule
WFG1  & $86.82{\:\pm\:0.68}$(4)  & $85.50{\:\pm\:0.24}$(5)  & $96.83{\:\pm\:0.24}$(3)  & $97.12{\:\pm\:0.22}$(2)  & $\textbf{97.91}{\:\pm\:}0.56$(1) \\
WFG2  & $109.60{\:\pm\:0.26}$(2)  & $109.38{\:\pm\:0.18}$(3)  & $\textbf{114.13}{\:\pm\:}3.76$(1)  & $108.79{\:\pm\:7.79}$(4)  & $100.13{\:\pm\:3.64}$(5) \\
WFG3  & $115.50{\:\pm\:0.27}$(2)  & $115.48{\:\pm\:0.17}$(3)  & $\textbf{116.42}{\:\pm\:}0.01$(1)  & $114.77{\:\pm\:0.34}$(4)  & $112.88{\:\pm\:0.67}$(5) \\
WFG4  & $\textbf{110.95}{\:\pm\:}0.26$(1)  & $109.17{\:\pm\:0.48}$(2)  & $105.99{\:\pm\:1.66}$(5)  & $107.09{\:\pm\:0.74}$(3)  & $106.07{\:\pm\:1.84}$(4) \\
WFG5  & $108.42{\:\pm\:0.99}$(3)  & $103.45{\:\pm\:1.21}$(5)  & $\textbf{110.33}{\:\pm\:}0.97$(1)  & $109.65{\:\pm\:1.29}$(2)  & $105.73{\:\pm\:1.55}$(4) \\
WFG6  & $113.15{\:\pm\:0.25}$(2)  & $109.64{\:\pm\:0.73}$(3)  & $\textbf{114.28}{\:\pm\:}0.04$(1)  & $109.58{\:\pm\:2.00}$(4)  & $109.49{\:\pm\:2.93}$(5) \\
WFG7  & $112.93{\:\pm\:0.48}$(4)  & $113.06{\:\pm\:0.36}$(3)  & $\textbf{114.33}{\:\pm\:}0.03$(1)  & $113.19{\:\pm\:0.59}$(2)  & $112.41{\:\pm\:0.51}$(5) \\
WFG8  & $109.72{\:\pm\:0.29}$(2)  & $109.27{\:\pm\:0.26}$(3)  & $\textbf{111.22}{\:\pm\:}0.22$(1)  & $109.17{\:\pm\:1.29}$(4)  & $105.98{\:\pm\:2.27}$(5) \\
WFG9  & $\textbf{110.70}{\:\pm\:}1.71$(1)  & $108.58{\:\pm\:0.57}$(3)  & $109.27{\:\pm\:0.66}$(2)  & $106.95{\:\pm\:2.02}$(4)  & $101.46{\:\pm\:2.87}$(5) \\
\midrule Rank & 2.33 (2)& 3.33 (4)& \textbf{1.78 (1)}& 3.22 (3)& 4.33 (5)\\
\bottomrule

\end{tabular}


\end{adjustbox}
\vspace*{-0.5cm}
\end{table}


\section{Discussion}
\label{sec:discussion}
We performed gradient-based multi-objective (MO) optimization by formulating the problem as a high-dimensional single-objective optimization problem based on the uncrowded hypervolume (UHV). We presented how the gradient of the UHV can be computed from the gradients of the MO function using the chain rule. 
We further showed that UHV gradient optimization can be solved with existing gradient ascent schemes, obtaining results competitive to or better than EAs and another gradient-based algorithm that performs hypervolume optimization. 
Future studies should additionally compare the presented UHV gradient-based algorithms to popular dominance-based EAs for HV optimization and investigate scalability also in $n$, the dimensionality of the underlying MO problem.

We have shown that the UHV is an effective and efficient approach for obtaining a set of non-dominated MO-solutions, that requires little to no extra care during the optimization process. 
In \cite{wang2017steering}, different techniques for steering dominated points are compared, and the uncrowded distance (UD) we use here is reminiscent of dominated point handling techniques such as secant slope weighting or gap-filling. However, a diversity loss is noted there as a possible disadvantage of gap-filling over the domination ranking technique employed in HIGA-MO, but we did not observe this in our results (Figure~\ref{fig:exp1}). 

UHV gradient ascent is not very sensitive to the initial step size, scales better when $p$ is large (Figure~\ref{fig:exp_maxhv_best}), and achieves a better hypervolume than EA-based UHV-optimization while requiring significantly fewer function evaluations (Figure~ \ref{tab:ppsn20_grad_wfg}).
When gradient information of the MO problem is missing, finite difference gradient approximation can be used, requiring $n+1$ MO-evaluations. Even with such approximations, UHV gradient ascent is competitive or even outperforms EAs on smaller computational budgets (Table~\ref{tab:ppsn20_grad_wfg}), although ultimately EAs often outperform the gradient-based algorithms. The effect of approximation errors in the finite difference approximation in UHV gradient ascent is negligible when using the Adam gradient scheme \cite{kingma2014adam}, as it was developed for stochastic gradient descent in which exact gradients are unavailable (or too expensive to compute). 

Any algorithm based on the hypervolume is limited by the hypervolume's computational complexity increasing in the number of objectives $m > 2$, e.g., $O(p^{m-2}\log(p))$ \cite{fonseca2006improved}. However, the approximation sets we consider are rather small (e.g., $p = 9$), and UHV gradient ascent shows good scalability in $p$, which encourages investigating cases with $m > 2$ objectives. 

Finally, as any gradient-based algorithm, UHV gradient ascent suffers from the risk of ending up in local optima.  A future research direction therefore is to hybridize UHV gradient ascent and EAs. Hybridization is however not trivial \cite{schutze2016hypervolume,bosman2011gradients} as both EAs and gradient-based algorithms rely on information of preceding iterations, and interleaving different algorithms might disrupt these mechanisms.


\vspace*{-0.5cm}
\subsubsection*{Acknowledgments.}
\small
This work was supported by the Dutch Research Council (NWO) through Gravitation Programme Networks 024.002.003 and is part of the research programme Open Technology Programme with project number 15586, which is financed by NWO, Elekta, and Xomnia. Further, this project is co-funded by the public-private partnership allowance for top consortia for knowledge and innovation (TKIs) from the Ministry of Economic Affairs.

\bibliographystyle{splncs04}
\bibliography{ppsn_bib.bib}

\end{document}